# Semi-parametric estimation of shifts

**Fabrice Gamboa and Jean-Michel Loubes**

*Laboratoire de Statistique et de Probabilités*
*UMR C5583, Université Paul Sabatier*
*118 Route de Narbonne, F-31062 Toulouse Cedex*
**e-mail:**
Fabrice.Gamboa@math.ups-tlse.fr; Jean-Michel.Loubes@math.ups-tlse.fr

**Elie Maza**

*Laboratoire Symbiose et Pathologie de Plantes*
*Pôle de Biotechnologie Végétale, INP-ENSAT*
*18 Chemin de Borde Rouge, BP 32607 Auzeville-Tolosane*
*31326 Castanet-Tolosan Cedex - France*
**e-mail:** Elie.Maza@ensat.fr

**Abstract**We observe a large number of functions differing from each other only by a translation parameter. While the main pattern is unknown, we propose to estimate the shift parameters using $M$-estimators. Fourier transform enables to transform this statistical problem into a semi-parametric framework. We study the convergence of the estimator and provide its asymptotic behavior. Moreover, we use the method in the applied case of velocity curve forecasting.

## 1. Introduction

A main issue in data mining is the feature extraction of a large set of curves. Indeed, classification methods enable to split the data into different homogeneous groups, each representing a specific mass behavior. But, within one group, the observations differ slightly the one from another. Such variations take into account the variability of the individuals inside one group. More precisely, there is a mean pattern such that, each observation curve is warped from this archetype by a warping function, see for examples [18].

In this work, we focus on the particular case where the individuals usually experience similar events, which are explained by a common pattern, but the starting time of the event occurs sooner or later. Classification methods, like repeated measures ANOVA or Principal Components Analysis of curves, see for instance [17], ignore this type of variability. Hence, computing a representative curve, for each group, severely distorts the analysis of the data. Indeed, the average curve (usually the mean or the median) oversmooths the studied phenomenon, and is not a good description of reality.

In our work, we restrict ourselves to the case where all the curves can be deduced the one from another by a shift parameter. Hence, we consider the following model: for $j = 1, \ldots, J$ and $i = 1, \ldots, n_j$, we observe

$$Y_{ij} = f\left(t_{ij} - \theta_j^*\right) + \sigma \varepsilon_{ij}, \tag{1}$$

where, $J$ stands for the number of curves $\mathcal{C}_j$, while $n_j$ is the number of observations for the $j$-th individual. Values $t_{ij}$ are observation times, which are assumed to be known. The unobserved warping effects $\theta_j^*$, $j = 1, \ldots, J$, are shift parameters which translate the unknown function $f$. We also choose $t_{ij} = t_i$ and $n_j = n$, which means that all curves are observed at the same time with the same occurrence. The errors $\varepsilon_{ij}$ for $(i,j) \in \{1, \ldots, n\} \times \{1, \ldots, J\}$ are i.i.d. with distribution $\mathcal{N}(0,1)$. Moreover, without loss of generality, we assume in the following that $\sigma = 1$ (see Remark 3.2). We aim at estimating the shift parameters $\theta_j^*$, $j = 1, \ldots, J$, in order to find a good representative of the feature $f$.

A more general problem has been tackled in the literature and some work has been done to find a representative of a large sample of close enough functions $f_j$, $j = 1, \ldots, J$, see for examples [18]. Indeed, in a general case, we observe realizations $y_{ij}$, $j = 1, \ldots, J$, $i = 1, \ldots, n_j$, from model

$$Y_{ij} = f_j(t_{ij}) + \varepsilon_{ij}, \tag{2}$$

where, $\varepsilon_{ij}$, $j = 1, \ldots, J$, $i = 1, \ldots, n_j$, are i.i.d. random variables, representing the observation noise. Hence, such functions $f_j$, $j = 1, \ldots, J$, are close from each other in the sense that there exists an unknown archetype $f$ and unknown warping functions $h_j$, $j = 1, \ldots, J$, such that, for all $j = 1, \ldots, J$,

$$\forall t \in [0, T],\ f_j(t) = f \circ h_j(t).$$

Examples of such data might be growth curves, longitudinal data in medicine, speech signals, traffic data or expenditure curves for some goods in the econometric domain. Our main motivation in this paper is the analysis of the vehicle speed evolution on a motorway. The data are curves, describing the evolution, on observation cells, of the daily vehicle speed. After performing classification procedures (see for instance [14] for a complete study), we obtain clusters of functions, each one representing a typical common behavior. Indeed, all the curves can be deduced one from another by a shift parameter.

This kind of issue led several statisticians to apply transformations to functions in order to get rid of the shifts and to align the curves. If a parametric model would be available a priori, the analysis would be made easier. But, if the data are numerous, there is not generally enough knowledge to build such a model. Thus, they turn into a non parametric framework. When the pattern is known, the problem turns to align a noisy observation with a fixed feature. Piccioni, Scarlatti and Trouvé in [15], Kneip, Li, MacGibbon and Ramsay in [11], or Ramsay and Li in [16] proposed curve registration methods. Their main idea is to align each curve on a target curve $f_0$, which means finding, for all $j \in \{1, \ldots, J\}$, the warping function $h_j$ minimizing

$$F_\lambda(f_0, f_j; h_j) = \int \|f_j \circ h_j(t) - f_0(t)\|^2 \mathrm{d}t + \lambda \int w_j^2(t) \mathrm{d}t,$$

where $h_j$ belongs to a particular smooth monotone family defined by the solution of the differential equation $D^2 h_j = w_j D h_j$. Hence, $w_j$ is simply $D^2 h_j / D h_j$, the relative curvature of $h_j$. Thus, penalizing $w_j$ yields both smoothness and monotonicity of $h_j$ (see [16] for more details). The main drawback of such methods is that they assume that the archetype $f_0$ is known, which is a reasonable assumption in pattern recognition, but which is unrealistic when the observed phenomenon is not well known as in our study. Alternatively, in a non parametric point of view, the pattern is replaced by its estimate. In this case, the issue is a matter of synchronizing sample curves. Wang and Gasser in [21] use kernel estimators. In another work, Gasser and Kneip, in [5], align the curves by aligning the local extrema of the functions, which are estimated as zeroes of the non parametric estimate of the derivative. In all cases, the issue of estimating the shifts is blurred by the estimation of the curves, which leads to non parametric rates of convergence.

Hence, it seems natural to study our regression problem (1) in a semi-parametric framework: the shifts are the parameters of interest to be estimated, while the pattern stands for an unknown nuisance functional parameter. A very general semi-parametric regression model called *Self-modelling regression (SEMOR)* has been considered in [10]. The model is $f_j(\cdot) = f(\cdot, \theta_j^*)$, $j \in \{1, \ldots, J\}$, and a general backfitting algorithm is studied. Roughly speaking, after initializing an estimate of $f$ by a first guess (using for example a kernel method), this algorithm is based on two recursive steps. In the first step, the estimation of $\theta_j^*$, $j = 1, \ldots, J$, is performed. In the second step, the estimate of $f$ is updated. In both steps, estimations are performed using a least squares criterion. In [10] a complete study, including the asymptotic normality of the estimates, is performed for the *Shape-invariant model (SIM)* introduced in [12]. See also [13], [8] and [9] for related works. Actually, the model studied in our paper (regression model (1)) falls in the SIM frame, so that, the methods studied in [10] may be applied. Nevertheless, the estimation procedure developed here is new, structurally simpler and computationaly easier to implement than the complicated backfitting algorithms.

The difficulty of the work is that the estimation must not rely on the pattern, even if the quantities are deeply linked. That is the reason why we will use an $M$-estimator built on the Fourier series of the data. Under identifiability assumptions, we provide a consistent method (Theorem 2.1) to estimate at the parametric rate of convergence the shifts $\theta_j^*$, $j = 1, \ldots, J$, when $f$ is unknown, and we show that fluctuations of the estimates are asymptotically Gaussian (Theorem 3.1). Further, our estimation method leads to a fast algorithm to align shifted curves without any prior assumption on the feature, due to semi-parametric techniques. We point out that this study can be linked first with the study of Golubev in [7], dealing with the semi-parametric efficiency in the estimation of shifts in a continuous observation scheme, and also with the study of Gassiat and Lévy-Leduc in [6], dealing with the estimation of the periodicity of a signal. Further, the mixed effects model (1) with random shifts is studied in [1] (see also [2]). We outline the fact that the method we propose handles a large variety of curves with minimal smoothness properties, namely we only require $L^1$ conditions.

The present paper falls into six parts. Section 2 is devoted to the definition of the model and to the description of the estimation method. In Section 3, we provide asymptotic properties of the estimators. As a matter of fact, we show that the estimators are convergent and asymptotically Gaussian. The estimating method is effectively performed in Section 4, on some simulated data, and then used to analyze road traffic data. We compare our results to another existing method. The technical lemmas and the proofs are gathered in Section 5 and Section 6.

## 2. Semi-parametric estimation of shifts

### 2.1. Model

For the $j$-th curve ($j = 1, \ldots, J$), we get $n$ observations $y_{ij}$, $i = 1, \ldots, n$, measured at equispaced times $t_i = \frac{i-1}{n}T \in [0, T[$, with $T \in \mathbb{R}_+^*$. We model these observations in the following way:

$$Y_{ij} = f(t_i - \theta_j^*) + \varepsilon_{ij}, \ j = 1, \ldots, J, \ i = 1, \ldots, n, \qquad (3)$$

where, $f : \mathbb{R} \to \mathbb{R}$ is an unknown $T$-periodic function, $\theta^* = (\theta_1^*, \ldots, \theta_J^*) \in \mathbb{R}^J$ is an unknown shift parameter, $\theta_j^*$ is the shift of the $j$-th curve, and, $\varepsilon_{ij}$, $i = 1, \ldots, n$, is a Gaussian white noise, with variance 1. For sake of simplicity, we consider an unitary variance, but all our results are still valid for a general variance.

Our aim is to estimate the translation factors $(\theta_j^*)$ without the knowledge of the pattern $f$. Due to the special structure of the model, Fourier analysis is well suited to conduct such a study, since the Fourier basis diagonalizes any translation. Then, using a Discrete Fourier Transform we may transform the model (3) into the following one (supposing $n$ is odd):

$$d_{jl} = e^{-il\alpha_j^*} c_l(f) + w_{jl}, \ j = 1, \ldots, J, \ l = -(n-1)/2, \ldots, (n-1)/2, \qquad (4)$$

where, $c_l(f) = \frac{1}{n} \sum_{m=1}^{n} f(t_m) e^{-i2\pi \frac{ml}{n}}$, $l = -(n-1)/2, \ldots, (n-1)/2$, are the discrete Fourier coefficients and $\alpha_j^* = \frac{2\pi}{T} \theta_j^* \in \mathbb{R}$, $j = 1, \ldots, J$, are the phase factors, and, for all $j \in \{1, \ldots, J\}$, $w_{jl}$, $l = -(n-1)/2, \ldots, (n-1)/2$, is a complex Gaussian white noise, with complex variance $1/n$, and with independent real and imaginary parts. As previously, our goal is to estimate the phase factors $\alpha_j^*$, $j = 1, \ldots, J$, without the knowledge of the Fourier coefficients of function $f$. Stricto sensu, the discrete Fourier coefficients are not the Fourier coefficients of the functions, but the bias induced is similar to the bias induced by any discretization in regression, which vanishes under some regularity assumptions, as shown in [4]. Hence, from now, we will consider the model (4) with $c_l(f) = \frac{1}{T} \int_0^T f(t) e^{-i2\pi \frac{tl}{T}} dt$. Observe that in this last equation we only have to assume that $f$ is integrable. Hence, if we only consider the discretized version given in Model 2.2, only a minimal smoothness conditions ($f \in L^1(\mathbb{R})$) is necessary, contrary to other methods in statistics, which require stronger conditions.

We point out that we are facing a semi-parametric model. As a matter of fact, we aim at estimating the parameter $\alpha^* = (\alpha_1^*, \cdots, \alpha_J^*)$ which depends on an unknown nuisance functional parameter $(c_l(f))_{l \in \mathbb{Z}}$, the Fourier coefficients of the unknown function $f$.

## 2.2. Identifiability

We notice that the model (4) is not identifiable for all translation parameters. Indeed, replacing $\alpha^*$ by

$$\begin{pmatrix} \alpha_1 \\ \vdots \\ \alpha_J \end{pmatrix} = \begin{pmatrix} \alpha_1^* \\ \vdots \\ \alpha_J^* \end{pmatrix} + c \begin{pmatrix} 1 \\ \vdots \\ 1 \end{pmatrix} + 2\pi \begin{pmatrix} k_1 \\ \vdots \\ k_J \end{pmatrix}, \, c \in \mathbb{R}, \, \begin{pmatrix} k_1 \\ \vdots \\ k_J \end{pmatrix} \in \mathbb{Z}^J, \quad (5)$$

and replacing $f(\cdot)$ by $f(\cdot - c)$, let invariant the equation (4). So, in order to ensure identifiability of the model, we restrict the parameter space $A$:

$$\begin{array}{ll} \text{i)} & A \text{ is compact,} \\ \text{ii)} & \alpha^* \in A, \\ \text{iii)} & \text{if } \alpha \in A \text{ and (5) holds for } \alpha, \text{ then } \alpha = \alpha^*. \end{array} \quad (6)$$

In this paper, we will mainly consider, in the Theorem 3.1, the parameter set $A_1 = \{\alpha \in [-\pi, \pi[^J : \alpha_1 = 0\}$. Hence, in (5) the constant $c$ must be equal to 0. Our fluctuations theorem can easily be transposed to other choices of parameter spaces, for example $A_2 = \left\{\alpha \in [-\pi, \pi[^J : \sum_{j=1}^J \alpha_j = 0 \text{ and } \alpha_1 \in \left[0, \frac{2\pi}{J}\right[\right\}$. In this last case, the condition $\sum_{j=1}^J \alpha_j = 0$ implies in (5) that $c = -\frac{2\pi}{J} \sum_{j=1}^J k_j$. So that, with equation (5), we can write that

$$\begin{pmatrix} \alpha_1 \\ \vdots \\ \alpha_J \end{pmatrix} = \begin{pmatrix} \alpha_1^* \\ \vdots \\ \alpha_J^* \end{pmatrix} + \frac{2\pi}{J} \begin{pmatrix} (J-1) & -1 & \cdots & -1 \\ -1 & \ddots & \ddots & \vdots \\ \vdots & \ddots & \ddots & -1 \\ -1 & \cdots & -1 & (J-1) \end{pmatrix} \begin{pmatrix} k_1 \\ \vdots \\ k_J \end{pmatrix}.$$

Hence, we get $J$ different solutions in $[-\pi, \pi[^J \subset \mathbb{R}^J$, and a unique solution with the additional condition $\alpha_1 \in \left[0, \frac{2\pi}{J}\right[$.

## 2.3. Estimation

Since we want to estimate the shifts without prior knowledge of the function $f$, we will consider a semi-parametric method, relying on an $M$-estimation procedure. Hence, the functional parameter is a nuisance parameter that does not play a role in the rate of converge of the estimates of the parameters, regardless of smoothness conditions for $f$.

For this, define, for any $\alpha = (\alpha_1, \ldots, \alpha_J) \in A$, the rephased coefficients

$$\tilde{c}_{jl}(\alpha) = e^{il\alpha_j} d_{jl}, \; j = 1, \ldots, J, \; l = -(n-1)/2, \ldots, (n-1)/2,$$

and the mean of these rephased Fourier coefficients

$$\hat{c}_l(\alpha) = \frac{1}{J} \sum_{j=1}^{J} \tilde{c}_{jl}(\alpha), \; l = -(n-1)/2, \ldots, (n-1)/2.$$

We have that $\tilde{c}_{jl}(\alpha^*) = c_l(f) + e^{il\alpha_j^*} w_{jl}$, for all $j \in \{1, \ldots, J\}$, and

$$\hat{c}_l(\alpha^*) = c_l(f) + \frac{1}{J} \sum_{j=1}^{J} e^{il\alpha_j^*} w_{jl}.$$

Hence, $|\tilde{c}_{jl}(\alpha) - \hat{c}_l(\alpha)|^2$ should be small when $\alpha$ is close to $\alpha^*$.

Now, consider a bounded measure $\mu$ on $[0, T]$ and set

$$\delta_l := \int_{[0,T]} \exp\left(\frac{2i\pi l}{T}\omega\right) \mathrm{d}\mu(\omega) \; (l \in \mathbb{Z}).$$

Obviously, the sequence $(\delta_l)$ is bounded. Without loss of generality we will assume that $\delta_0 = 0$. Assume further that $\sum_l |\delta_l|^2 |c_l(f)|^2 < +\infty$. So that $f * \mu$ is a well defined square integrable function:

$$f * \mu(x) = \int f(x-y) \mathrm{d}\mu(y).$$

Consider the following empirical contrast function:

$$M_n(\alpha) = \frac{1}{J} \sum_{j=1}^{J} \sum_{l=-\frac{n-1}{2}}^{\frac{n-1}{2}} |\delta_l|^2 \, |\tilde{c}_{jl}(\alpha) - \hat{c}_l(\alpha)|^2 . \tag{7}$$

In the sequel, we will always assume that:

The set $\{l : \delta_l c_l \neq 0\}$ contains at least two integers which are coprime. (8)

The random function $M_n$ is non negative. Furthermore, its minimum value should be reached close to the true parameter $\alpha^*$. Then, the following theorem provides the consistency of the $M$-estimator, defined by

$$\hat{\alpha}_n = \arg\min_{\alpha \in A} M_n(\alpha).$$

**Theorem 2.1** *Under the following assumptions on $f$ and on the weight sequence $(\delta_l)_{l\in\mathbb{Z}}$:*

$$\begin{cases} \sum_{l\in\mathbb{Z}} |\delta_l|^2 |c_l(f)|^2 & < +\infty, \\ \sum_{l\in\mathbb{Z}} |\delta_l|^4 |c_l(f)|^2 & < +\infty \end{cases} \quad (9)$$

*we have that $\hat{\alpha}_n \xrightarrow[n\to+\infty]{\mathbf{P}_{\alpha^*}} \alpha^*$.*

We point out that we only assume that $f*\mu$ and $f*\mu*\mu$ are square integrable and yet are able to build estimates of the shifts in the model (4). The computation of the estimator is quick since only a Fast Fourier Transform algorithm and a minimization algorithm of a quadratic functional are needed.

**Proof 2.2 (Proof of Theorem 2.1)** *The proof of this theorem follows the classical guidelines of the convergence of $M$-estimators (see for example [19]). Indeed, the contrast is split into two parts, a determinist and a random one. Then, it suffices to show that the following conditions hold for the criterion function to ensure consistency of $\hat{\alpha}_n$.*

*i) Convergence to a contrast function:*

$$M_n(\alpha) \xrightarrow[n\to+\infty]{\mathbf{P}_{\alpha^*}} K(\alpha),\ \alpha \in A, \quad (10)$$

*where $K(\cdot)$ has a unique minimum at $\alpha^*$.*

*ii) Set the modulus of uniformly continuity $W$, defined by*

$$W(n,\eta) = \sup_{\|\alpha-\beta\|\leq\eta} |M_n(\alpha) - M_n(\beta)|.$$

*There exists two sequences $(\eta_k)_{k\in\mathbb{N}}$ and $(\epsilon_k)_{k\in\mathbb{N}}$, decreasing to zero, such that for a large enough $k$, we have*

$$\lim_{n\to+\infty} \mathbf{P}_{\alpha^*}\left(W(n,\eta_k) > \epsilon_k\right) = 0. \quad (11)$$

*These two conditions are fulfilled, as it is proved in Section 6. Notice that we chose to privilege the uniform convergence of the modulus of continuity of the contrast and not the uniform convergence of the criterion itself. Nevertheless, the two proofs use the same kind of arguments, i.e proving the uniform convergence of empirical processes of Gaussian variables.*

## 3. Asymptotic normality

In this section, we prove that the estimator built in the previous section is asymptotically Gaussian, and we give its asymptotic covariance matrix. In general, the asymptotic covariance matrix hardly depends on the geometric structure of the parameter space $A$. So, for sake of simplicity, we study the asymptotic normality for the parameter space $A_1$. Hence, the parameter space has dimension $J-1$, and we rewrite this set as $\tilde{A}_1 = [-\pi, \pi[^{J-1}$ and, any element in $\tilde{A}_1$ as $\alpha = (\alpha_2, \ldots, \alpha_J)$. Also, for sake of simplicity, in this section and in the proofs of Theorem 3.1, we will write $M_n(\alpha)$ instead of $M_n(0, \alpha_2, \ldots, \alpha_J)$. So, we consider any estimator defined by

$$\hat{\alpha}_n = \arg \min_{\alpha \in \tilde{A}_1} M_n(\alpha).$$

**Theorem 3.1** *Under the following assumptions* $(\delta_l)_{l \in \mathbb{Z}}$:

$$\begin{cases} 0 < \sum_{l \in \mathbb{Z}} |\delta_l|^2 l^2 |c_l(f)|^2 & < \infty \\ \sum_{l \in \mathbb{Z}} |\delta_l|^2 l^4 |c_l(f)|^2 & < \infty \\ \sum_{l=-n}^{n} |\delta_l|^4 l^4 & = o(n^2), \end{cases} \quad (12)$$

*we get that*

$$\sqrt{n}(\hat{\alpha}_n - \alpha^*) \xrightarrow[n \to +\infty]{\mathcal{D}} \mathcal{N}_{J-1}(0, \Gamma), \quad (13)$$

*with*

$$\Gamma = \frac{\sum_{l \in \mathbb{Z}} |\delta_l|^4 l^2 |c_l(f)|^2}{\left(\sum_{l \in \mathbb{Z}} |\delta_l|^2 l^2 |c_l(f)|^2\right)^2} (I_{J-1} + U_{J-1}),$$

*where, $I_{J-1}$ is the identity matrix of dimension $J-1$, and $U_{J-1}$ is the square matrix of dimension $J-1$ whose all entries are equal to one.*

**Remark 3.2** *If the white noise in the model (3) has a variance equal to $\sigma^2$, then the limit distribution in the previous theorem has a covariance matrix equal to $\sigma^2 \Gamma$.*

**Proof 3.3 (Proof of Theorem 3.1)** *Recall that the $M$-estimator is defined as the minimum of the criterion function $M_n(\alpha)$. Hence, we get*

$$\nabla M_n(\hat{\alpha}_n) = 0,$$

*where $\nabla$ is the gradient operator. A second order decomposition leads to: there exists $\bar{\alpha}_n$ in a neighborhood of $\alpha^*$ such that*

$$\sqrt{n}(\hat{\alpha}_n - \alpha^*) = -\left[\nabla^2 M_n(\bar{\alpha}_n)\right]^{-1} \sqrt{n} \nabla M_n(\alpha^*), \quad (14)$$

*where $\nabla^2$ is the Hessian operator. Now, using the two asymptotic results from Proposition 5.1 and from Proposition 5.2, we get*

$$\sqrt{n} \nabla M_n(\alpha^*) \xrightarrow[n \to +\infty]{\mathcal{D}} \mathcal{N}_{J-1}(0, \Gamma_0),$$

$$\left[\nabla^2 M_n(\bar{\alpha}_n)\right]^{-1} \xrightarrow[n \to +\infty]{\mathbf{P}_{\alpha^*}} V,$$

where $V$ is a non negative symmetric matrix of dimension $J - 1$. Hence, if we set $\Gamma = V'\Gamma_0 V$, the result of Theorem 3.1 follows easily. Finally, we see that $\sum_{l=-n}^{n} |\delta_l|^4 l^4 = o(n^2)$ implies that $\sum_{l=-n}^{n} |\delta_l|^4 l^2 = o(n)$. Indeed, $\frac{1}{n^2} \sum_{l=-n}^{n} |\delta_l|^4 l^4 \geq n^2 |\delta_n|^4 + n^2 |\delta_{-n}|^4$, so $\lim_{|n| \to +\infty} |\delta_n|^2 |n| = 0$. Hence $\frac{1}{n} \sum_{l=-n}^{n} |\delta_l|^4 l^2 = o(1)$. Moreover $\sum_l |\delta_l|^2 l^4 |c_l(f)|^2 < +\infty$ implies Assumption (19). So that, the set of assumptions (12) implies the ones of both Propositions 5.2 and 5.1.

Observe that the extra terms $(\delta_l)_{l \in \mathbb{Z}}$ used in the definition (7) smooth the criterion function $M_n(\alpha)$. Indeed, without this term, i.e under the choice $\delta_l = 1$, the random part of the derivative of the criterion function does not converge towards a determinist function but to a random process, which prevents the study of the asymptotic distribution. The weights enable to get rid of this part, smoothing the contrast to zero.

Moreover, the convergence of the criterion is speeded up by using these smoothing weights. We illustrate this purpose on Section 4 by comparing a weighted criterion with a non weighted one (see Figure 3). Moreover, this result will be highlighted in the proof of Theorem 3.1.

**Practical choice of the $\delta_l$'s**

The problem of choosing the weights $(\delta_l)_{l \in \mathbb{Z}}$ in the definition of the criterion function is important. If we work with $L^2$ functions, the assumption (12) is satisfied for example as soon as $|\delta_l| = O\left(|l|^{-2-\nu}\right)$, for some $\nu > 0$. In the simulations our functions are much more regular. Hence, we have taken $\delta_l = 1/|l|^{1.3}$. This choice guarantees consistency and good numerical results. Moreover, in order to illustrate the importance of the weight sequence, we have also taken $\delta_l \equiv 1$ and $\delta_l = 1/|l|^2$ in Figure 3. Indeed, when looking at the asymptotic variance, we can see that there is a trade-off which leads to a lower bound for the smoothing sequence, the smaller the weights, the larger the variance. Since the function $f$ is unknown and so the sequence does not depend on the Fourier coefficients, hence the optimal choice for $(\delta_l)_{l \in \mathbb{Z}}$ should be given by semi-parametric efficiency. Using Cauchy-Schwarz's inequality, we get that

$$\frac{\sum_{l \in \mathbb{Z}} |\delta_l|^4 l^2 |c_l(f)|^2}{\left(\sum_{l \in \mathbb{Z}} |\delta_l|^2 l^2 |c_l(f)|^2\right)^2} \geq \left(\sum_{l \in \mathbb{Z}} l^2 |c_l(f)|^2\right)^{-1}.$$

This case, corresponding to the least favorable case in the semi-parametric efficiency framework, is obtained for the optimal choice of coefficients $\delta_l = 1$. If an asymptotic fluctuation results would hold for this sequence, we would obtain:

$$\sqrt{n} \sqrt{\sum_{l \in \mathbb{Z}} l^2 |c_l(f)|^2} \left(\hat{\alpha}_n - \alpha^*\right) \xrightarrow[n \to +\infty]{\mathcal{D}} \mathcal{N}_{J-1}\left(0, I_{J-1} + U_{J-1}\right).$$

Nevertheless, for the choice of the weight sequence $\delta_l = 1$ the asymptotic normality result does not hold. Non optimality as regards asymptotic efficiency is the price to pay both to deal with a discretized version of the regression model and to handle simultaneous estimation for all the unknown functions.

Maybe, a different way of estimation could get rid of this drawback. Yet, another choice could have been done to smooth the contrast by restricting the number of Fourier coefficients, as it is done in [7] for example. Some links could also be established between the estimator we consider and a Bayesian penalized maximum likelihood estimator, where the weights $(\delta_l)_{l \in \mathbb{Z}}$ stand for a particular choice of a prior over the unknown function $f$. This Bayesian point of view is tackled in [3]. However, the optimal choice of the smoothing parameter to obtain efficiency is a very interesting issue in the semi-parametric framework (i.e when the weights are not allowed to depend on the Fourier coefficients of the functions). Quite posterior to the first submission of this work, this problem has been solved in a [20].

**Remark 3.4** *Throughout all the work, we assume that the observation noise in the model* (3) *is Gaussian. Nevertheless, we could get rid of this assumption with moment conditions on the errors.*

## 4. Applications and simulations

In this section, we present some numerical applications of the method. The first one gives results on simulated data. The second one is based on an experiment on human fingers force. The last one is carried out with traffic data.

The optimization algorithm used in any resolution is based on a Krylov method (the conjugate gradient method). Indeed, minimizing an $L^2$ criterion function with a conjugate gradient algorithm yields a reduced step number, and hence, a small complexity.

**Simulated data**

Simulated data are carried out as follows:

$$y_{ij} = f(t_i - \theta_j^*) + \varepsilon_{ij}, \; j = 1, \ldots, J, \; i = 1, \ldots, n,$$

with the following choice of parameters: $J = 10$; $n = 100$; values $t_i = -\pi + \frac{i-1}{n} 2\pi$, $i = 1, \ldots, n$, are equally spaced points on $[-\pi, \pi[$; $f(t) = 15 \sin(4t)/(4t)$; $(\theta_2^*, \ldots, \theta_J^*)$ are simulated with a uniform law on $[-\pi/4, \pi/4]$ and $\theta_1^* = 0$; for all $j \in \{1, \ldots, J\}$, for all $i \in \{1, \ldots, n\}$, values $\varepsilon_{ij}$ are simulated from a Gaussian law with mean 0 and standard deviation 1. Results are given on Figure 1. The target function $f$ is considered as a $2\pi$-periodic function ($T = 2\pi$), hence $\alpha^* = \theta^*$. The function $f$ is plotted by a solid line in Figure 1 (d). Figure 1 (a) shows simulated data $y_{ij}$, $j = 1, \ldots, J$, $i = 1, \ldots, n$. The cross-sectional mean curve of these data is given on Figure 1 (d) by the dotted line. We can see that this mean function is representativeness of data. Indeed, the amplitude of higher optimum is reduced, and smallest ones disappeared. Figure 1 (c) shows curves unshifted by the estimated parameters. The mean function of these unshifted curves is given on Figure 1 (d) by dashed line. Figure 1 (b) plots $\theta_j^*$ on abscissa axis against $\hat{\theta}_j$ on ordinate axis, $j = 1, \ldots, J$. Estimations are very close to true parameters. Comparison between mean curves, before and after the shift estimation, is straightforward.

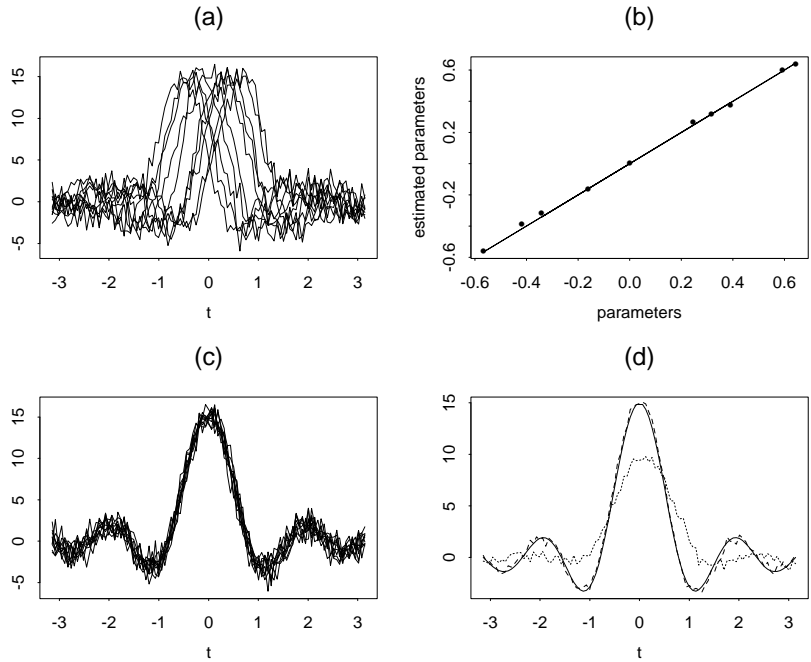

Figure 1. *Estimation results with the M-estimation methodology.*

We now compare our estimations with those obtained with an existing method: curve registration by landmarks. This method aims at aligning curves by, first, estimating landmarks of curves (here, the maximum) and by, secondly, aligning these landmarks. For more details on this procedure, see [5]. In Figure 2, we show the results on our simulated data. These results are not as good as those we obtain with our method. It can be explained by the fact that we need first to estimate each curve maximum by a non parametric method (a kernel estimation), which leads to estimation errors. Moreover, our method uses all information given by the data, not only that given by landmarks.

In order to illustrate the importance of the weight sequence, we compare now the criterion function for various values of $(\delta_l)_{l \in \mathbb{Z}}$. For this purpose, simulated data sets are carried out, with $J = 2$, $\theta_1^* = 0$ and $\theta_2^* = \pi/3$. Figure 3 shows the obtained results. The first column of this figure presents these simulated data sets, with respectively, $\sigma = 1$ in Figure 3 (a,1), $\sigma = 3$ in Figure 3 (a,2), $\sigma = 5$ in Figure 3 (a,3) and $\sigma = 7$ in Figure 3 (a,4). The second column presents the unweighted criterion functions, i.e with $\delta_l \equiv 1$, associated respectively with (a,1), (a,2), (a,3) and (a,4). The third and fourth columns present the associated weighted criterion functions $M_n(\cdot)$ with, respectively, $\delta_l = 1/|l|^{1.3}$ and $\delta_l = 1/|l|^2$. In all these figures, the vertical dashed line represents the value $\pi/3$ where the minimum value of our criterion function is achieved. It clearly appears that without the weight sequence, i.e with $\delta_l \equiv 1$, the criterion function converges to a random process. Moreover, the variance of this random process is proportional to the noise variance $\sigma^2$. We also see that even with an important noise variance, as in Figure 3 (a,4), our weighted criterion functions, (c,4) and (d,4), are smooth, with a unique minimum around $\pi/3$. This shows that our procedure is quite robust to the SNR. Moreover, it appears that the impact of the exponent $\beta > 1.25$ of the weight sequence $\delta_l = 1/|l|^\beta$ is only on the amplitude of the $M$-function.

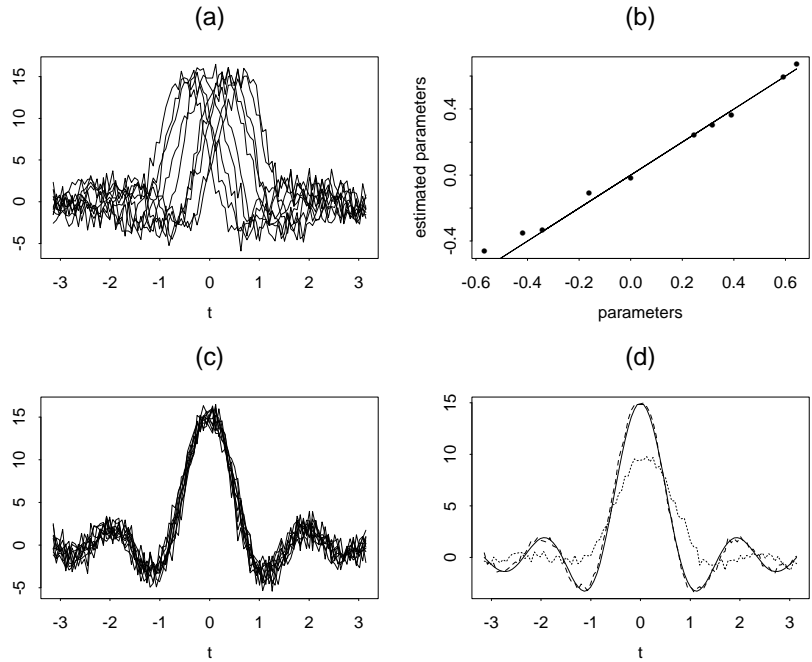

Figure 2. *Estimation results with the landmark methodology.*

These numerical results emphasize the fact that the weight sequence $(\delta_l)_{l \in \mathbb{Z}}$ is important, but that its value can be easily chosen.

### Pinch force data

Data presented here are extracted from an experiment described in [16] with a Curve Registration methodology. Data represent the force exerted by the thumb and forefinger on a force meter during 20 brief pinches. These 20 force measurements having arbitrary beginning, Ramsay and Li in [16] begin their study by a landmark alignment of curve maxima (with single shifts). These aligned data are shown in Figure 4 (a).

Our purpose is to study these data with our shift estimation methodology. Shift estimations and unshifted curves are respectively shown in Figure 4 (b) and

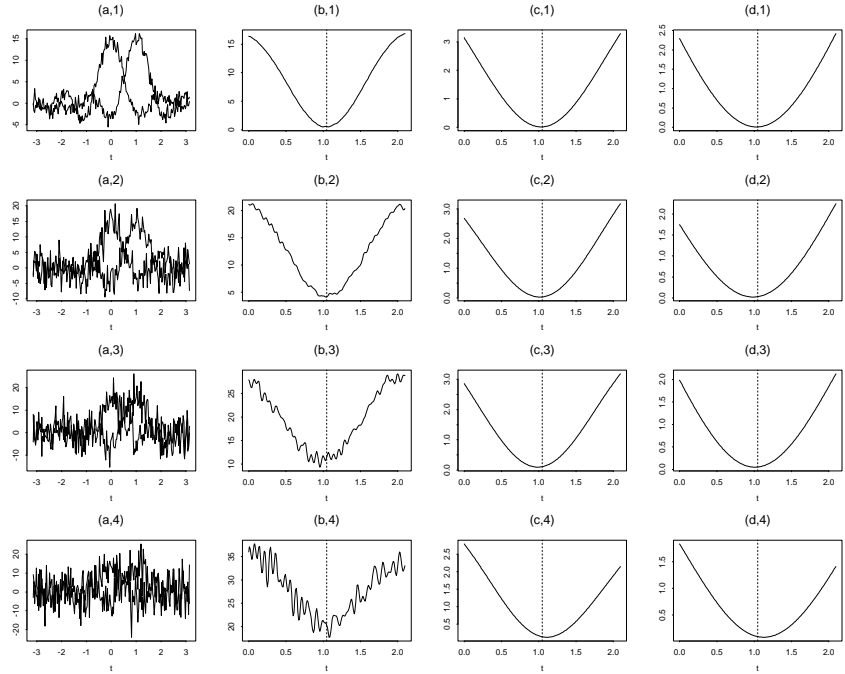

FIGURE 3. *Criterion functions, with different values of the weight sequence* $(\delta_l)_{l\in\mathbb{Z}}$.

Figure 4 (c). In Figure 4 (b), we only show a boxplot of the estimated parameters because, obviously, we do not know the real parameters. We note that shift parameters are almost all close to zero, between $-10^{-3}$ and $3 \times 10^{-3}$. In this case, landmark alignment unshift quite well the data. Nevertheless, comparing Figure 4 (a) and Figure 4 (c), we can see that curves are slightly better aligned after our shift estimation methodology. In Figure 4 (d), the cross-sectional mean curves of unshifted curves (solid line) and of primary curves (dotted line) are almost the same ones.

### Application to road traffic forecasting

Most of the Parisian road traffic network is equipped with a traffic road measurement infrastructure. The main elements of this infrastructure are counting stations. These sensors are situated approximately every 500 meters on main trunk roads (motorways and speedways principally). Every counting station measures, daily, the average speed of vehicle flow on 6 minutes periods. We consider measurements from 5 AM to 11 PM, hence, the length of the daily measurement is 180. We note $y_{ij}$ the speed measurement of day $j \in \{1, \ldots, J\}$ and of period $i \in \{1, \ldots, n\}$, with $n = 180$.

Our purpose is to improve, with our shift estimation methodology, an existing forecasting methodology. This forecasting methodology is described in [14]. This procedure is based on a classification method. We dispose of a sample of $J$ speed curves and we want to summarize it by a small number $N$ of standard profiles, representatives of each cluster.

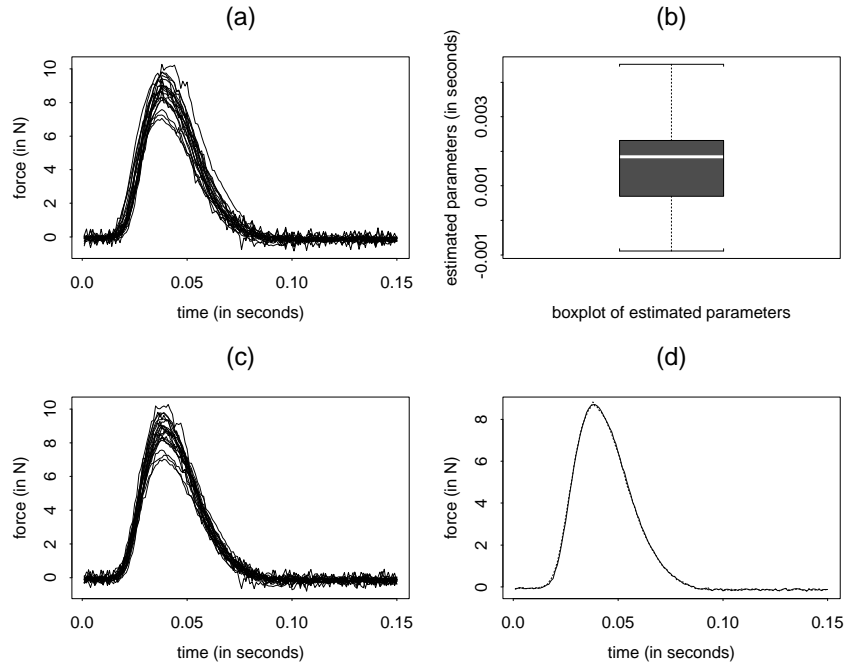

Figure 4. *Shift estimation results on the pinch force data set.*

Consider several clusters of $J$ speed curves. Indeed, we note frequently that many subgroups are composed by curves describing the same behavior. For example, we observe a speed curve subgroup with a same traffic jam or speed reduction, but with different starting times for each curve. Thus, Figure 5 (a) represent a particular cluster on a particular counting station. Figure 5 (b) is a boxplot of the estimated shifts $\hat{\theta}_j$, $j = 1, \ldots, J$. Unshifted curves are plotted on Figure 5 (c). So, in this homogeneous cluster, where only a shift phenomenon appear, the mean curves in Figure 5 (d) of unshifted curves (solid line) and of primary curves (dotted line) aren't the same. The shift estimated mean is clearly more representative of individual pattern.

## 5. Technical Lemmas

The two following propositions, Proposition 5.1 and Proposition 5.2, are used in the proof of asymptotic normality (Theorem 3.1). Their proofs are postponed to the appendix.

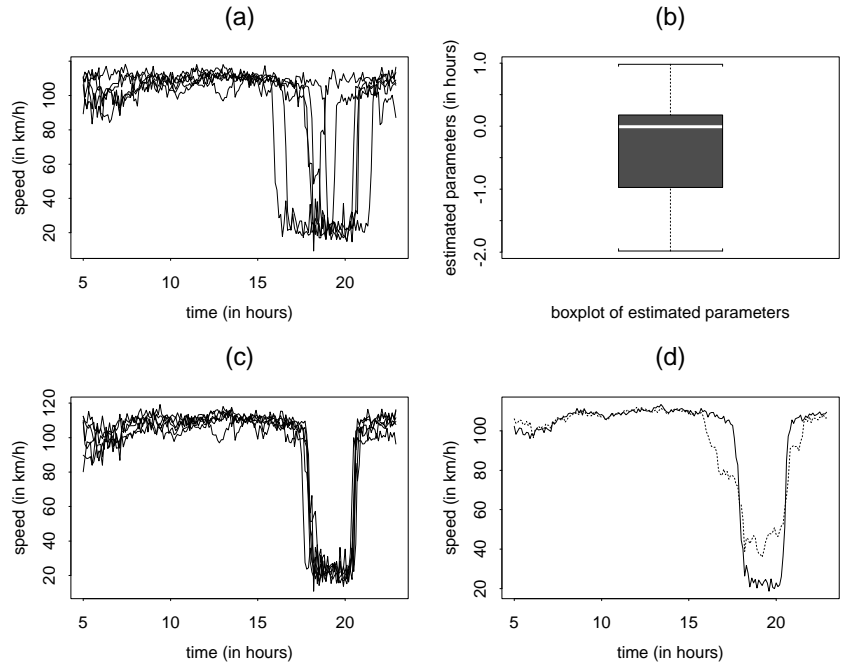

FIGURE 5. *Shift estimation results on a particular traffic data set.*

**Proposition 5.1** *Assume that the $\delta_l$'s are such that*

$$\sum_{l \in \mathbb{Z}} l^2 |\delta_l|^4 |c_l(f)|^2 < +\infty \tag{15}$$

$$\sum_{l \in \mathbb{Z}} l^2 |\delta_l|^4 = o(n) \tag{16}$$

*Then*

$$\sqrt{n} \nabla M_n(\alpha^*) \xrightarrow[n \to +\infty]{\mathcal{D}} \mathcal{N}_{J-1}(0, \Gamma), \tag{17}$$

*where the variance matrix is $\Gamma = \frac{4}{J^2} \sum_{l \in \mathbb{Z}} |\delta_l|^4 l^2 |c_l(f)|^2 \left( I_{J-1} - \frac{1}{J} U_{J-1} \right).$*

**Proposition 5.2** *Assume moreover that the sequence $\delta_l$ satisfies*

$$\sum_{l \in \mathbb{Z}} l^2 |\delta_l|^2 |c_l(f)|^2 < +\infty \tag{18}$$

$$\sum_{l \in \mathbb{Z}} l^4 |\delta_l|^4 |c_l(f)|^2 < +\infty \tag{19}$$

$$\sum_{l \in \mathbb{Z}} l^4 |\delta_l|^4 = o(n^2) \tag{20}$$

Then, for any sequence $(\bar{\alpha}_n)_{n\in\mathbb{N}}$ such that $\|\bar{\alpha}_n - \alpha^*\| \leq \|\hat{\alpha}_n - \alpha^*\|$ $(n \in \mathcal{N})$, we have

$$\nabla^2 M_n(\bar{\alpha}_n) \xrightarrow[n\to+\infty]{\mathbf{P}_{\alpha^*}} \frac{2}{J^2} \sum_{l\in\mathbb{Z}} |\delta_l|^2 l^2 |c_l(f)|^2 (JI_{J-1} - U_{J-1}). \tag{21}$$

## 6. Appendix

Let $z$ be a complex number and $\bar{z}$ its conjugate. We write $\mathfrak{Re}(z) = \frac{1}{2}(z+\bar{z})$ (the real part of $z$) and $\mathfrak{Im}(z) = \frac{1}{2i}(z-\bar{z})$ (the imaginary part of $z$).

**Proof 6.1 (Proof of Theorem 2.1)** *In the sequel we assume without loss of generality that all the $\alpha_j^*$ are equal to 0. Consider the following notation: for all $j = 1, \ldots, J$, for all $l = -(n-1)/2, \ldots, (n-1)/2$, $w_{jl} = \frac{1}{\sqrt{n}} \xi_{jl} = \frac{1}{\sqrt{n}} \left(\xi_{jl}^x + i\xi_{jl}^y\right)$. Here, $\left(\xi_{jl}^x\right)$ and $\left(\xi_{jl}^y\right)$ are independent Gaussian sequences, with law $\mathcal{N}_n(0, I_n)$. Also, set*

$$\forall l = -(n-1)/2, \ldots, (n-1)/2, \ c_l(f) = |c_l(f)|e^{i\theta_l}, \ \text{with } \theta_l \in [0, 2\pi[.$$

*Note that*

$$M_n(\alpha) = \sum_{l\in\mathbb{Z}} \delta_l^2 \left[ \frac{1}{J} \sum_{j=1}^{J} |\tilde{c}_{jl}(\alpha)|^2 - |\hat{c}_l(\alpha)|^2 \right].$$

*Using the following decompositions*

$$|\tilde{c}_{jl}(\alpha)|^2 = |c_l(f)|^2 + |w_{jl}|^2 + 2\mathfrak{Re}[c_l(f)\bar{w}_{jl}],$$

$$|\hat{c}_l(\alpha)|^2 = \left|\frac{1}{J}\sum_{j=1}^{J} e^{il\alpha_j} c_l(f)\right|^2 + \left|\frac{1}{J}\sum_{j=1}^{J} e^{il\alpha_j} w_{jl}\right|^2$$

$$+ 2\mathfrak{Re}\left[\left(\frac{1}{J}\sum_{j=1}^{J} e^{il\alpha_j} c_l(f)\right)\left(\frac{1}{J}\sum_{k=1}^{J} e^{-il\alpha_k} \bar{w}_{kl}\right)\right],$$

*lead to the following expression of the criterion function $M_n(\alpha)$*

$$M_n(\alpha) = \sum_{l=-\frac{n-1}{2}}^{\frac{n-1}{2}} |\delta_l|^2 |c_l(f)|^2 - \sum_{l=-\frac{n-1}{2}}^{\frac{n-1}{2}} \left|\frac{1}{J}\sum_{j=1}^{J} e^{il(\alpha_j)} \delta_l c_l(f)\right|^2 \tag{22}$$

$$+ \frac{J-1}{J^2} \sum_{j=1}^{J} \frac{1}{n} \sum_{l=-\frac{n-1}{2}}^{\frac{n-1}{2}} |\delta_l|^2 \left(\xi_{jl}^{x\,2} + \xi_{jl}^{y\,2}\right) \tag{23}$$

$$- \frac{2}{J^2} \sum_{j=1}^{J} \sum_{k>j} \frac{1}{n} \sum_{l=-\frac{n-1}{2}}^{\frac{n-1}{2}} |\delta_l|^2 \Big[\cos(l(\alpha_j - \alpha_k))(\xi_{jl}^x \xi_{kl}^x + \xi_{jl}^y \xi_{kl}^y)$$

$$+ \sin(l(\alpha_j - \alpha_k))(\xi_{jl}^x \xi_{kl}^y - \xi_{jl}^y \xi_{kl}^x)\Big] \tag{24}$$

$$+ \frac{2(J-1)}{J^2} \sum_{j=1}^{J} \frac{1}{\sqrt{n}} \sum_{l=-\frac{n-1}{2}}^{\frac{n-1}{2}} |\delta_l|^2 |c_l(f)| \left[\cos(\theta_l)\xi_{jl}^x + \sin(\theta_l)\xi_{jl}^y\right]$$

$$- \frac{2}{J^2} \sum_{j=1}^{J} \sum_{k \neq j} \frac{1}{\sqrt{n}} \sum_{l=-\frac{n-1}{2}}^{\frac{n-1}{2}} |\delta_l|^2 |c_l(f)| \left[\cos(l(\alpha_j - \alpha_k) + \theta_l)\xi_{kl}^x \right.$$
$$\left. + \sin(l(\alpha_j - \alpha_k) + \theta_l)\xi_{kl}^y\right]. \tag{25}$$

We have split the criterion function into four different terms: (22), (23), (24), and (25). We aim at proving the convergence of these terms to a determinist contrast function and the uniform convergence of their increments.

- The term (22) is a deterministic one. Using Parseval theorem, we have that

$$(22) \xrightarrow[n \to +\infty]{} \int_0^T |(f * \mu)(t)|^2 \frac{dt}{T} - \int_0^T \left| \frac{1}{J} \sum_{j=1}^{J} (f * \mu)(t + \theta_j - \theta_j^*) \right|^2 \frac{dt}{T}.$$

- The term (23) is a pure noise term composed of terms of the type $\frac{1}{n} \sum_{l=-\frac{n-1}{2}}^{\frac{n-1}{2}} |\delta_l|^2 \xi_{jl}^{x\,2}$. Since $\xi_{jl}$, $l \in \mathcal{I}_n$ are independent, by the SLLN we get

$$\frac{1}{n} \sum_{l=-\frac{n-1}{2}}^{\frac{n-1}{2}} |\delta_l^2| \xi_{jl}^{x\,2} \xrightarrow[n \to +\infty]{\mathbf{P}_{\alpha^*}} \Lambda(\delta) < +\infty,$$

for a constant $\Lambda(\delta)$ which only depends on the choice of the smoothing sequence and not on the unknown parameter $\alpha$. This constant is bounded since the weights are bounded. Note that if $\mu$ has a density lying in $\mathbb{L}^2$, this constant vanishes since $1/n \sum_{l \in \mathbb{Z}} |\delta_l|^2 \to 0$.

- The first term in (24) is also a pure noise term composed of terms of the type

$$U_n(\alpha_j - \alpha_k) = \frac{1}{n} \sum_{l=-\frac{n-1}{2}}^{\frac{n-1}{2}} |\delta_l|^2 \cos(l(\alpha_j - \alpha_k))\xi_{jl}^x \xi_{kl}^y, \ k \neq j.$$

One has $\mathbb{E}(U_n) = 0$ and $\mathbb{E}(\exp(t \sum_{l=-\frac{n-1}{2}}^{\frac{n-1}{2}} \delta_l^2 |\xi_{jl}^x \xi_{kl}^y|)) < +\infty$. Thus for $a, b > 0$ a Bernstein type inequality holds giving

$$\mathbf{P}_{\alpha^*} \left( \sup_{\alpha_j - \alpha_k \in [-2\pi, 2\pi] \cap \mathbb{Z}/n^2} |U_n(\alpha_j - \alpha_k)| > x \right) \leq O\left(n^2 \exp(-\frac{1}{2} \frac{n^2 x^2}{a + bnx})\right).$$

Choosing $x_n = \sqrt{8b \log(n)/n}$, we get both

$$x_n \to 0 \quad \text{and} \quad \mathbf{P}_{\alpha^*}(\sup_{\alpha_j - \alpha_k \in [-2\pi, 2\pi] \cap \mathbb{Z}/n^2} |U_n(\alpha_j - \alpha_k)| > x_n) \to 0.$$

*For $|h| \leq 1/n^2$ the inequality $|\cos lh - 1| \leq 1/n$ leads to*

$$|U_n(\alpha_j - \alpha_k + h) - U_n(\alpha_j - \alpha_k)| \leq \frac{1}{n^2} \sum_{l=-\frac{n-1}{2}}^{\frac{n-1}{2}} |\delta_l|^2 \xi^x_{jl} \xi^y_{kl}$$

$$\leq \frac{1}{n} \left( \frac{1}{n} \sum_{l=-\frac{n-1}{2}}^{\frac{n-1}{2}} \xi^{x\,2}_{jl} \right)^{\frac{1}{2}} \left( \frac{1}{n} \sum_{l=-\frac{n-1}{2}}^{\frac{n-1}{2}} |\delta_l|^2 \xi^{x\,2}_{kl} \right)^{\frac{1}{2}} \to 0,$$

*so that*

$$\sup_{\alpha_j - \alpha_k \in [-2\pi, 2\pi]} |U_n(\alpha_j - \alpha_k)| \xrightarrow[n \to +\infty]{\mathbf{P}_{\alpha^*}} 0.$$

*Hence the first term in (24) goes to 0 in probability.*

- *The remaining term in (24) is similar to the term (25), which has the same asymptotic behavior as*

$$V_n(\alpha_j - \alpha_k) = \frac{1}{\sqrt{n}} \sum_{l=-\frac{n-1}{2}}^{\frac{n-1}{2}} |\delta_l|^2 |c_l(f)| \cos(l(\alpha_j - \alpha^*_j - \alpha_k) + \theta_l) \xi^x_{kl}.$$

*Here we only give the proof of the uniform convergence of (25) which holds under slight modifications of the second term of (24). As the noise is Gaussian we get*

$$\begin{cases} \mathbb{E}(V_n(\alpha_j - \alpha_k)) = 0 \\ \mathrm{Var} V_n(\alpha_j - \alpha_k) \leq \frac{2}{n} \sum_{l=-\frac{n-1}{2}}^{\frac{n-1}{2}} |\delta_l|^4 |c_l(f)|^2 \leq \frac{2}{n} \sum_{l=-\infty}^{+\infty} |\delta_l|^4 |c_l(f)|^2 \end{cases}$$

*so that*

$$\mathbf{P}_{\alpha^*} \left( \sup_{\alpha_j - \alpha_k \in [-2\pi, 2\pi]} |V_n(\alpha_j - \alpha_k)| \geq \sqrt{\frac{16 \log n \sum_{l=-\infty}^{+\infty} |\delta_l|^2 |c_l(f)|^2}{n}} \right)$$

$$\leq K n^2 \exp(-4 \log n).$$

*Using again previous bound for $h \leq 1/n^2$, we obtain*

$$|V_n(\alpha_j - \alpha_k + h) - V_n(\alpha_j - \alpha_k)| \leq \frac{1}{n\sqrt{n}} \sum_{l=-\frac{n-1}{2}}^{\frac{n-1}{2}} |\delta_l| |c_l(f)| \xi^x_{kl}$$

$$\leq \frac{1}{n} \left( \sum_{l=-\infty}^{+\infty} |\delta_l|^2 |c_l(f)|^2 \right)^{\frac{1}{2}} \left( \frac{1}{n} \sum_{l=-\infty}^{+\infty} \xi^{x\,2}_{kl} \right)^{\frac{1}{2}}$$

$$\to 0$$

*so that*

$$\sup_{\alpha_j - \alpha_k \in [-2\pi, 2\pi]} |V_n(\alpha_j - \alpha_k)| \xrightarrow[n \to +\infty]{\mathbf{P}_{\alpha^*}} 0.$$

*In conclusion, we have that* $\sup_{\alpha \in [-\pi, \pi]^J} |M_n(\alpha) - K(\alpha) - | \xrightarrow[n \to +\infty]{\mathbf{P}_{\alpha^*}} 0$ *with*

$$K(\alpha) = \int_0^T |(f * \mu)(t)|^2 \frac{dt}{T} - \int_0^T \left| \frac{1}{J} \sum_{j=1}^J (f * \mu)(t + \theta_j - \theta_j^*) \right|^2 \frac{dt}{T} + \Lambda(\delta).$$

*This ensures that (11) is fulfilled and that the convergence property in (10) is ensured. It remains to be seen that the asymptotic contrast enables to identify the $\alpha_j$'s, which concludes the proof of Condition (10). Cauchy-Schwartz inequality yields that*

$$\int_0^T \left| \frac{1}{J} \sum_{j=1}^J (f * \mu)(t + \theta_j - \theta_j^*) \right|^2 \frac{dt}{T}$$

$$\leq \int_0^T \frac{1}{J} \sum_{j=1}^J |(f * \mu)(t + \theta_j - \theta_j^*)|^2 \frac{dt}{T} = \int_0^T |(f * \mu)(t)|^2 \frac{dt}{T},$$

*hence, $K(\cdot) \geq 0$, and the minimum value is reached for*

$$\int_0^T |(f * \mu)(t)|^2 \frac{dt}{T} = \int_0^T \left| \frac{1}{J} \sum_{j=1}^J (f * \mu)(t + \theta_j - \theta_j^*) \right|^2 \frac{dt}{T},$$

*which is equivalent, using Parseval theorem to*

$$\sum_{l \in \mathbb{Z}} |\delta_l c_l(f)|^2 = \sum_{l \in \mathbb{Z}} \left| \frac{1}{J} \sum_{j=1}^J \delta_l c_l(f) e^{il(\alpha_j - \alpha_j^*)} \right|^2. \tag{26}$$

*So, we have that*

$$(26) \iff \forall l \in \{l : c_l \delta_l \neq 0\}, \left| \frac{1}{J} \sum_{j=1}^J e^{il(\alpha_j - \alpha_j^*)} \right|^2 = 1.$$

*Now, from (8) this implies that*

$$\forall j = 1, \ldots, J, \ \alpha_j = \alpha_j^* + c \, [2\pi], \ c \in \mathbb{R}.$$

*In a matrix way, we get the equation (5), i.e*

$$\begin{pmatrix} \alpha_1 \\ \vdots \\ \alpha_J \end{pmatrix} = \begin{pmatrix} \alpha_1^* \\ \vdots \\ \alpha_J^* \end{pmatrix} + c \begin{pmatrix} 1 \\ \vdots \\ 1 \end{pmatrix} + 2\pi \begin{pmatrix} k_1 \\ \vdots \\ k_J \end{pmatrix}, \ c \in \mathbb{R}, \ \begin{pmatrix} k_1 \\ \vdots \\ k_J \end{pmatrix} \in \mathbb{Z}^J.$$

Hence, since $\alpha \in A$ and $A$ is defined by (6), we have shown that $\alpha_j = \alpha_j^*$ for all $j = 1, \ldots, J$. Since $\alpha \mapsto K(\alpha)$ achieves its unique minimum for $\alpha = \alpha^*$, the condition (10) is fulfilled.

**Proof 6.2 (Proof of Proposition 5.1)** *The first and the second derivatives of the empirical contrast, for all $k \in \{2, \ldots, J\}$, for all $m \in \{2, \ldots, J\}$ can be written as:*

$$\frac{\partial M_n}{\partial \alpha_k}(\alpha) = \frac{2}{J} \sum_{l=-\frac{n-1}{2}}^{\frac{n-1}{2}} |\delta_l|^2 l \mathfrak{Im}\left(\tilde{c}_{kl}(\alpha)\overline{\hat{c}_l(\alpha)}\right), \qquad (27)$$

$$\frac{\partial^2 M_n}{\partial \alpha_k^2}(\alpha) = \frac{2}{J^2} \sum_{l=-\frac{n-1}{2}}^{\frac{n-1}{2}} |\delta_l|^2 l^2 \mathfrak{Re}\left(\tilde{c}_{kl}(\alpha) \sum_{j \neq k} \overline{\tilde{c}_{jl}(\alpha)}\right), \qquad (28)$$

$$\forall m \neq k, \; \frac{\partial^2 M_n}{\partial \alpha_k \partial \alpha_m}(\alpha) = -\frac{2}{J^2} \sum_{l=-\frac{n-1}{2}}^{\frac{n-1}{2}} |\delta_l|^2 l^2 \mathfrak{Re}\left(\tilde{c}_{kl}(\alpha)\overline{\tilde{c}_{ml}(\alpha)}\right). \qquad (29)$$

*By straightforward calculations, we get that*

$$\sqrt{n}\frac{\partial M_n}{\partial \alpha_k}(\alpha^*) = \frac{2}{J} \sum_{l=-\frac{n-1}{2}}^{\frac{n-1}{2}} |\delta_l|^2 l \left(|c_l(f)|\left(V_l^k - V_l\right) + W_l^k\right),$$

*where, for all $l \in \mathbb{Z}$,*

$$W_l^k = \frac{1}{J\sqrt{n}} \sum_{j=1}^{J} \Big[\sin(l(\alpha_k^* - \alpha_j^*))(\xi_{kl}^x \xi_{jl}^x + \xi_{kl}^y \xi_{jl}^y)$$
$$+ \cos(l(\alpha_k^* - \alpha_j^*))(\xi_{kl}^y \xi_{jl}^x - \xi_{kl}^x \xi_{jl}^y)\Big],$$

$$V_l^k = (\cos(l\alpha_k^* + \theta_l)\xi_{kl}^x - \sin(l\alpha_k^* + \theta_l)\xi_{kl}^y), \; \text{and } V_l = \frac{1}{J}\sum_{j=1}^{J} V_l^j.$$

*Let, for $l \in \mathbb{Z}$, $Y_l = (\xi_{1l}^x \xi_{2l}^x \cdots \xi_{Jl}^x \xi_{1l}^y \xi_{2l}^y \cdots \xi_{Jl}^y)'$, and, let $f_l^k$ be the vector of length $2J$, defined by $(f_l^k)_k = \cos(l\alpha_k^* + \theta_l)$, $(f_l^k)_{J+k} = -\sin(l\alpha_k^* + \theta_l)$, and $(f_l^k)_i = 0$ if $i \notin \{k, J+k\}$. As a consequence, we get the following expression for $V_l^k$: $V_l^k = \langle f_l^k, Y_l \rangle = f_l^{k'} Y_l$. In a same way, for $l \in \mathbb{Z}$, let $\bar{B}_l^k$ be the $(2J) \times (2J)$ matrix defined by rows by*

$$\left(\bar{B}_l^k\right)_k = (\sin[l(\alpha_k^* - \alpha_1^*)] \cdots \sin[l(\alpha_k^* - \alpha_J^*)] \; -\cos[l(\alpha_k^* - \alpha_1^*)] \cdots$$
$$- \cos[l(\alpha_k^* - \alpha_J^*)]),$$

$$\left(\bar{B}_l^k\right)_{J+k} = (\cos[l(\alpha_k^* - \alpha_1^*)] \cdots \cos[l(\alpha_k^* - \alpha_J^*)] \; \sin[l(\alpha_k^* - \alpha_1^*)] \cdots \sin[l(\alpha_k^* - \alpha_J^*)]),$$

$$\left(\bar{B}_l^k\right)_i = (0 \cdots 0) \; \text{if } i \notin \{k, J+k\}.$$

Further, let the symmetric matrix $B_l^k$ be defined by $B_l^k = \frac{\bar{B}_l^k + (\bar{B}_l^k)'}{2}$. Hence, write $W_l^k = \frac{1}{J\sqrt{n}} Y_l' B_l^k Y_l$. Define also for $k = 2, \ldots, J$ $\tilde{B}_l^k = \frac{2}{J} B_l^k$, $\tilde{f}_l^k = \frac{2}{J} f_l^k$.

Our aim is to study the asymptotic distribution of the gradient $\sqrt{n} \nabla M_n(\alpha^*)$. For this purpose consider $u = (u_2, \ldots, u_J)' \in \mathbb{R}^{J-1}$ and $t = (t_2, \ldots, t_J)' \in \mathbb{R}^{J-1}$, and define the couple of random variables:

$$(R_n, S_n) = \left( \frac{2}{J} \sum_{k=2}^J u_k \sum_{l=-\frac{n-1}{2}}^{\frac{n-1}{2}} |\delta_l|^2 l |c_l(f)| (V_l^k - V_l), \frac{2}{J} \sum_{k=2}^J t_k \sum_{l=-\frac{n-1}{2}}^{\frac{n-1}{2}} |\delta_l|^2 l W_l^k \right).$$

Using previous notations, we get

$$R_n = \sum_{l=-\frac{n-1}{2}}^{\frac{n-1}{2}} |\delta_l|^2 l |c_l(f)| \langle g_l(u), Y_l \rangle, \text{ with } g_l(u) = \sum_{k=2}^J u_k \left( \tilde{f}_l^k - \frac{1}{J} \sum_{j=1}^J \tilde{f}_l^j \right),$$

$$S_n = \frac{1}{J\sqrt{n}} \sum_{l=-\frac{n-1}{2}}^{\frac{n-1}{2}} |\delta_l|^2 l Y_l' A_l(t) Y_l, \text{ with } A_l(t) = \sum_{k=2}^J t_k \tilde{B}_l^k.$$

First note that $\mathbb{E}(S_n) = 0$. Moreover Assumption (16) implies that

$$\text{Var} S_n \leq \frac{1}{n} \sum_{l \in \mathbb{Z}^*} |\delta_l|^4 |l|^2 \to 0.$$

Hence, $S_n \xrightarrow[n \to +\infty]{\mathbf{P}_{\alpha^*}} 0$. So, the quadratic part is vanishing in probability when $n$ increases.

For the last term, we have that $\langle g_l(u), Y_l \rangle \sim \mathcal{N}\left(0, \|g_l(u)\|_2^2\right)$, where $\|g_l(u)\|_2^2 = \frac{4}{J^2} u' \left( I_{J-1} - \frac{1}{J} U_{J-1} \right) u$, with, $I_{J-1}$ the $J-1$ identity matrix, and $U_{J-1}$ the $J-1 \times J-1$ matrix which all entries are equal to 1. The independence of the $Y_l$'s yields that under Assumption (15)

$$R_n \xrightarrow[n \to +\infty]{\mathcal{D}} \mathcal{N}\left(0, \frac{4}{J^2} \sum_{l \in \mathbb{Z}} |\delta_l|^4 l^2 |c_l(f)|^2 u' \left( I_{J-1} - \frac{1}{J} U_{J-1} \right) u \right).$$

Finally we get that

$$\sqrt{n} \nabla M_n(\alpha^*) \xrightarrow[n \to +\infty]{\mathcal{D}} \mathcal{N}_{J-1}\left(0, \frac{4}{J^2} \sum_{l \in \mathbb{Z}} |\delta_l|^4 l^2 |c_l(f)|^2 \left( I_{J-1} - \frac{1}{J} U_{J-1} \right) \right).$$

**Proof 6.3 (Proof of Proposition 5.2)** *First, we pay attention to the non diagonal terms of the matrix of the second derivatives. For $m \neq k$, we get after some calculations:*

$$-\frac{J^2}{2}\frac{\partial^2 M_n}{\partial \alpha_k \partial \alpha_m}(\alpha) = \sum_{l=-\frac{n-1}{2}}^{\frac{n-1}{2}} |\delta_l|^2 l^2 \mathfrak{Re}\left(\tilde{c}_{kl}(\alpha)\overline{\tilde{c}_{ml}(\alpha)}\right)$$

$$= \sum_{l=-\frac{n-1}{2}}^{\frac{n-1}{2}} |\delta_l|^2 l^2 |c_l(f)|^2 \cos\left(l[\alpha_k - \alpha_k^* + \alpha_m^* - \alpha_m]\right) \quad (30)$$

$$+ \sum_{l=-\frac{n-1}{2}}^{\frac{n-1}{2}} |\delta_l|^2 l^2 |c_l(f)| \left(\cos[l(\alpha_k - \alpha_k^* - \alpha_m) + \theta_l] w_{ml}^x \right.$$
$$\left. + \sin[l(\alpha_k - \alpha_k^* - \alpha_m) + \theta_l] w_{ml}^y\right) \quad (31)$$

$$+ \sum_{l=-\frac{n-1}{2}}^{\frac{n-1}{2}} |\delta_l|^2 l^2 |c_l(f)| \left(\cos[l(\alpha_m - \alpha_m^* - \alpha_k) + \theta_l] w_{kl}^x \right.$$
$$\left. + \sin[l(\alpha_m - \alpha_m^* - \alpha_k) + \theta_l] w_{kl}^y\right) \quad (32)$$

$$+ \sum_{l=-\frac{n-1}{2}}^{\frac{n-1}{2}} |\delta_l|^2 l^2 \left[\cos(l[\alpha_k - \alpha_m])(w_{kl}^x w_{ml}^x + w_{kl}^y w_{ml}^y) \right.$$
$$\left. - \sin(l[\alpha_k - \alpha_m])(w_{kl}^y w_{ml}^x - w_{kl}^x w_{ml}^y)\right]. \quad (33)$$

*We now study the asymptotic behaviour of each term separately. Indeed, the second derivatives are taken at a point $\bar{\alpha}_n$ which converges to $\alpha^*$: $\bar{\alpha}_n$ is in the neighborhood of $\alpha^*$ with radius $\|\alpha^* - \hat{\alpha}_n\|$. Hence, we need conditions to claim uniform convergence of $\nabla^2 M_n(\cdot)$.*

*First note that the deterministic term $\sum_{l=-(n-1)/2}^{(n-1)/2} |\delta_l|^2 l^2 |c_l(f)|^2 \cos(l[\bar{\alpha}_k - \alpha_k^* + \alpha_m^* - \bar{\alpha}_m])$ converges towards $\sum_l |\delta_l|^2 l^2 |c_l(f)|^2$ as soon as $\sum_l |\delta_l|^2 l^4 |c_l(f)|^2 < +\infty$, as assumed in (12).*

*Now consider the random terms. Since for all $k \in \{1, \ldots, J\}$, the random variables $w_{kl}^x$ and $w_{kl}^y$ follow a Gaussian law $\mathcal{N}(0, 1/n)$, we consider the independent variables $\xi_{kl}^x$ and $\xi_{kl}^y$ such that $w_{kl}^x = \frac{1}{\sqrt{n}}\xi_{kl}^x$ and $w_{kl}^y = \frac{1}{\sqrt{n}}\xi_{kl}^y$. For the two second terms (31) and (32), we write*

$$(31) = \frac{1}{\sqrt{n}}\sum |\delta_l|^2 l^2 |c_l(f)| \left(\cos[l(\alpha_k - \alpha_k^* - \alpha_m) + \theta_l]\xi_{ml}^x \right.$$
$$\left. + \sin[l(\alpha_k - \alpha_k^* - \alpha_m) + \theta_l]\xi_{ml}^y\right),$$
$$(32) = \frac{1}{\sqrt{n}}\sum |\delta_l|^2 l^2 |c_l(f)| \left(\cos[l(\alpha_m - \alpha_m^* - \alpha_k) + \theta_l]\xi_{kl}^x \right.$$
$$\left. + \sin[l(\alpha_m - \alpha_m^* - \alpha_k) + \theta_l]\xi_{kl}^y\right),$$

*whose asymptotic behaviours are of the same nature. As in the proof of Proposition 5.1, Condition (19) leads to $\sup_{\alpha \in \tilde{A}_1} (31) \xrightarrow[n \to +\infty]{\mathbf{P}_{\alpha^*}} 0$, and $\sup_{\alpha \in \tilde{A}_1} (32) \xrightarrow[n \to +\infty]{\mathbf{P}_{\alpha^*}} 0$. Further, Assumption (20) implies that (33) converges in probability uniformly*

to 0. The diagonal terms can be written as follows:

$$\frac{J^2}{2}\frac{\partial^2 M_n}{\partial \alpha_k^2}(\alpha) = \sum_{l=-\frac{n-1}{2}}^{\frac{n-1}{2}} |\delta_l|^2 l^2 \mathfrak{Re}\left(\tilde{c}_{kl} \sum_{j\neq k} \overline{\tilde{c}_{jl}}\right)$$

$$= \sum_{l=-\frac{n-1}{2}}^{\frac{n-1}{2}} |\delta_l|^2 l^2 |c_l(f)|^2 \sum_{j\neq k} \cos\left(l[\alpha_k - \alpha_k^* + \alpha_j^* - \alpha_j]\right) \quad (34)$$

$$+ \sum_{l=-\frac{n-1}{2}}^{\frac{n-1}{2}} |\delta_l|^2 l^2 |c_l(f)| \sum_{j\neq k} \Big(\cos[l(\alpha_k - \alpha_k^* - \alpha_j) + \theta_l] w_{jl}^x$$
$$+ \sin[l(\alpha_k - \alpha_k^* - \alpha_j) + \theta_l] w_{jl}^y\Big) \quad (35)$$

$$+ \sum_{l=-\frac{n-1}{2}}^{\frac{n-1}{2}} |\delta_l|^2 l^2 |c_l(f)| \sum_{j\neq k} \Big(\cos[l(\alpha_j - \alpha_j^* - \alpha_k) + \theta_l] w_{kl}^x$$
$$+ \sin[l(\alpha_j - \alpha_j^* - \alpha_k) + \theta_l] w_{kl}^y\Big) \quad (36)$$

$$+ \sum_{l=-\frac{n-1}{2}}^{\frac{n-1}{2}} |\delta_l|^2 l^2 \sum_{j\neq k} \Big[\cos(l[\alpha_k - \alpha_j])(w_{kl}^x w_{jl}^x + w_{kl}^y w_{jl}^y)$$
$$- \sin(l[\alpha_k - \alpha_j])(w_{kl}^y w_{jl}^x - w_{kl}^x w_{jl}^y)\Big]. \quad (37)$$

Using similar arguments as for the previous terms, we can see that, under the same assumptions we get that all the terms (34), (35), (36) and (37) converges uniformly, and we get

$$\frac{\partial^2 M_n}{\partial \alpha_k^2}(\bar{\alpha}_n) \xrightarrow[n\to+\infty]{\mathbf{P}_{\alpha^*}} \frac{2(J-1)}{J^2} \sum_{l\in\mathbb{Z}} |\delta_l|^2 l^2 |c_l(f)|^2.$$

As a result, gathering the two previous results leads to the following asymptotic behavior:

$$\nabla^2 M_n(\bar{\alpha}_n) \xrightarrow[n\to+\infty]{\mathbf{P}_{\alpha^*}} \frac{2}{J^2} \sum_{l\in\mathbb{Z}} |\delta_l|^2 l^2 |c_l(f)|^2 \left(JI_{J-1} - U_{J-1}\right),$$

which proves the result. Moreover, this matrix is invertible. As a result, we have that

$$\left[\nabla^2 M_n(\bar{\alpha}_n)\right]^{-1} \xrightarrow[n\to+\infty]{\mathbf{P}_{\alpha^*}} \frac{J}{2\sum_{l\in\mathbb{Z}} |\delta_l|^2 l^2 |c_l(f)|^2} \left(I_{J-1} + U_{J-1}\right).$$

**Acknowledgement**

We warmly thank the anymous referee for his careful reading and helpful coments.